\documentclass[reqno]{amsart}
\usepackage{amsmath}
\usepackage{amssymb}
\usepackage{amsthm}
\begin{document}
%%%%%%%%%%%%%%%%%%%%%%%%%%%%%%%%%%%%%%%%%%%%%%%%%%%%%%%%%%%%%%%%%%%%%%%
%%%%%%%%%%%%%%%%%%%%%%%%%     Macros      %%%%%%%%%%%%%%%%%%%%%%%%%%%%%
%%%%%%%%%%%%%%%%%%%%%%%%%%%%%%%%%%%%%%%%%%%%%%%%%%%%%%%%%%%%%%%%%%%%%%%
\def\eq#1{{\rm(\ref{#1})}}
\theoremstyle{plain}% default
\newtheorem{thm}{Theorem}[section]
\newtheorem{lem}[thm]{Lemma}
\newtheorem{prop}[thm]{Proposition}
\newtheorem{cor}[thm]{Corollary}
\theoremstyle{definition}
\newtheorem{dfn}[thm]{Definition}
\newtheorem{rem}[thm]{Remark}
\newtheorem{ex}{Example}
\def\Ker{\mathop{\rm Ker}}
\def\Coker{\mathop{\rm Coker}}
\def\ind{\mathop{\rm ind}}
\def\Re{\mathop{\rm Re}}
\def\vol{\mathop{\rm vol}}
\def\SO{\mathbin{\rm SO}}
\def\Im{\mathop{\rm Im}}
\def\min{\mathop{\rm min}}
\def\Spec{\mathop{\rm Spec}\nolimits}
\def\Hol{{\textstyle\mathop{\rm Hol}}}
\def\ge{\geqslant}
\def\le{\leqslant}
\def\C{{\mathbin{\mathbb C}}}
\def\R{{\mathbin{\mathbb R}}}
\def\N{{\mathbin{\mathbb N}}}
\def\Z{{\mathbin{\mathbb Z}}}
\def\D{{\mathbin{\mathcal D}}}
\def\H{{\mathbin{\mathcal H}}}
\def\M{{\mathbin{\mathcal M}}}
\def\al{\alpha}
\def\be{\beta}
\def\ga{\gamma}
\def\de{\delta}
\def\ep{\epsilon}
\def\io{\iota}
\def\ka{\kappa}
\def\la{\lambda}
\def\ze{\zeta}
\def\th{\theta}
\def\vp{\varphi}
\def\si{\sigma}
\def\up{\upsilon}
\def\om{\omega}
\def\De{\Delta}
\def\Ga{\Gamma}
\def\Th{\Theta}
\def\La{\Lambda}
\def\Om{\Omega}
\def\ts{\textstyle}
\def\sst{\scriptscriptstyle}
\def\sm{\setminus}
\def\na{\nabla}
\def\pd{\partial}
\def\op{\oplus}
\def\ot{\otimes}
\def\bigop{\bigoplus}
\def\iy{\infty}
\def\ra{\rightarrow}
\def\longra{\longrightarrow}
\def\dashra{\dashrightarrow}
\def\t{\times}
\def\w{\wedge}
\def\d{{\rm d}}
\def\bs{\boldsymbol}
\def\ci{\circ}
\def\ti{\tilde}
\def\ov{\overline}
\def\md#1{\vert #1 \vert}
\def\nm#1{\Vert #1 \Vert}
\def\bmd#1{\big\vert #1 \big\vert}
\def\cnm#1#2{\Vert #1 \Vert_{C^{#2}}} %% C^k or C^{k,\al} norms
\def\lnm#1#2{\Vert #1 \Vert_{L^{#2}}} %% L^q norms
\def\bnm#1{\bigl\Vert #1 \bigr\Vert}
\def\bcnm#1#2{\bigl\Vert #1 \bigr\Vert_{C^{#2}}} %% big C^k etc norms.
\def\blnm#1#2{\bigl\Vert #1 \bigr\Vert_{L^{#2}}} %% big L^q norms
%%%%%%%%%%%%%%%%%%%%%%%%%%%%%%%%%%%%%%%%%%%%%%%%%%%%%%%%%%%%%%%%%%%%%%%
%%%%%%%%%%%%%%%%%%%%%     Text of paper    %%%%%%%%%%%%%%%%%%%%%%%%%%%%
%%%%%%%%%%%%%%%%%%%%%%%%%%%%%%%%%%%%%%%%%%%%%%%%%%%%%%%%%%%%%%%%%%%%%%%
\title[Asymptotically Cylindrical Ricci-Flat Manifolds]{Asymptotically Cylindrical Ricci-Flat Manifolds}
\author[Sema Salur]{Sema Salur$^*$}
\thanks{* Research supported in part by AWM-NSF Mentoring grant}
\address {Department of Mathematics, Northwestern University, IL, 60208 }
\email{salur@math.northwestern.edu }

\begin{abstract} Asymptotically cylindrical Ricci-flat
manifolds play a key role in constructing Topological Quantum
Field Theories. It is particularly important to understand their
behavior at the cylindrical ends and the natural restrictions on
the geometry. In this paper we show that an orientable, connected,
asymptotically cylindrical manifold $(M,g)$ with Ricci-flat metric
$g$ can have at most two cylindrical ends. In the case where there
are two such cylindrical ends then there is reduction in the
holonomy group Hol$(g)$ and $(M,g)$ is a cylinder.
\end{abstract}
\date{}
\maketitle
%%%%%%%%%%%%%%%%%%%%%%%%%%%%%%%%%%%%%%%%%%%%%%%%%%%%%%%%%%%%%%%%%%%%%
%%%%%%%%%%%%%%%%%%%%%%%%%%%%%%%%%%%%%%%%%%%%%%%%%%%%%%%%%%%%%%%%%%%%%
\section{Introduction}
\label{co1}

Let $(M^n,g)$ be a connected asymptotically cylindrical manifold
with Ricci-flat metric $g$. Examples of such manifolds are
asymptotically cylindrical Riemannian manifolds whose holonomy
groups Hol$(g)$ are subgroups of $SU(m)$, for $n=2m$, and
subgroups of the exceptional Lie group $G_2$, for $n=7$. Examples
of asymptotically cylindrical Riemannian 6-manifolds with holonomy
$SU(3)$ were first constructed by Tian and Yau \cite{TiYa1},
\cite{TiYa2} and Kovalev \cite{Kova}. In his paper, \cite{Kova},
Kovalev also showed that the Riemannian product of these
6-manifolds with a circle produced asymptotically cylindrical
Riemannian 7-manifolds with holonomy $SU(3)\subset G_2$. These are
the first examples of asymptotically cylindrical $G_2$ manifolds.

Other examples of noncompact Ricci-flat manifolds $(M^n,g)$ are
hyperk\"{a}hler manifolds, manifolds with holonomy group
Hol$(g)\subset Sp(m)$, for $n=4m$, and $Spin(7)$ manifolds with
holonomy group Hol$(g)\subset Spin(7)$, for $n=8$.

Asymptotically cylindrical Calabi-Yau and $G_2$ manifolds are the
main objects of the Topological Quantum Field Theories, as is
explained in \cite{Leun2}. In order to construct consistent
TQFT's, it is essential to understand the geometric restrictions
at the cylindrical ends of a given orientable, connected,
noncompact Ricci-flat manifold with some small decay rate. A very
natural and fundamental question is whether it is possible for
these manifolds to have multiple cylindrical ends.

In this paper we prove the following theorem.

\begin{thm} Let $(M,g)$ be an orientable, connected, asymptotically cylindrical
Riemannian manifold with Ricci-flat metric $g$ and have $l$
cylindrical ends. Then $l\leq2$, and in the case when $l=2$ there
is reduction in the holonomy group Hol$(g)$ and $(M,g)$ is a
cylinder. \label{co1thm}
\end{thm}

\begin{rem}

The proof of this theorem also follows from the Cheeger-Gromoll
splitting theorem \cite{CG}. The key ingredients of the splitting
theorem are the maximum principle for continuous functions,
Laplacian comparison theorems, constructions of rays, lines and
Busemann functions. For more details on the subject see
\cite{Petersen}. In this paper we give an alternative proof and
show that the reduction in holonomy can be obtained by just using
the analytic set-up for Fredholm properties of an elliptic
operator on noncompact manifolds. This analytic set-up was
developed by Lockhart and McOwen, \cite{LoMc} and by Melrose,
\cite{Mel1}, \cite{Mel2}.
\end{rem}

\section{Asymptotically cylindrical manifolds}
\label{co23} In this section we first introduce the definitions of
cylindrical and asymptotically cylindrical Riemannian manifolds.
\begin{dfn} An $n$-dimensional Riemannian manifold $(M_0,g_0)$ is called
{\it cylindrical} if $M_0=X\t\R$ and $g_0$ is compatible with this
product structure, that is,
\begin{equation*}
g_0=g_X+\d t^2,
\end{equation*}
where $X$ is a compact, connected $(n-1)$ manifold with Riemannian
metric $g_X$. \label{co2def1}
\end{dfn}

\begin{dfn} A connected, complete $n$-manifold $(M,g)$ with $l$
cylindrical ends is called {\it asymptotically cylindrical with
decay rate} ${\bf\be=(\be_1,\ldots ,\be_l)}\in\R^l$, $\be_j<0$ for
$1\leq j \leq l$, if there exist cylindrical $n$-manifolds
$(M_{0i},g_{0i})$ with $M_{0i}=X_i\t\R$ for connected $X_i$,
compact subsets $K_i\subset M$, a real number $R$, and
diffeomorphisms $\Psi_i:X_i\t(R,\iy)\ra M\sm K_i$, such that the
pull-back metric $(\Psi_i)^*(g)$ satisfies
$\bmd{\na^k_{0i}((\Psi_i)^*(g)-g_{0i})}=O(e^{\be_i t})$ on
$X_i\t(R,\iy)$ for all $k\ge 0$, where $\na_{0i}$ is the
Levi-Civita connection of the cylindrical metric $g_{0i}$ on
$M_{0i}=X_i\t\R$. \label{co2def22}
\end{dfn}

\begin{rem} For simplicity, we will assume that the decay rates
are all equal for each cylindrical end, that is $\be_j=\be_k$,
$1\leq j,k\leq l$.
\end{rem}

\subsection{An Example}
If there are additional geometric structures on the cylindrical
manifold $(M_0,g_0)=(X\t\R, g_0)$ then they should be also
compatible with this product structure. In addition to the metric
$g$, the geometric structures on the asymptotically cylindrical
manifold $(M,g)$ should converge to order $O(e^{\be t})$ to the
cylindrical geometric structures on $(M_0,g_0)$, with all of their
derivatives.

One example where we have such additional structures on $(M,g)$ is
an asymptotically cylindrical $G_2$ manifold. A 7-dimensional
Riemannian manifold $(M, g)$ is called a {\em $G_2$ manifold} if
the holonomy group of its Levi-Civita connection of $g$ lies
inside of $G_2\subset SO(7)$. A $G_2$ manifold is equipped with a
vector cross product $\times$ on its tangent bundle and a harmonic
(calibration) $3$-form $\varphi \in \Omega^{3}(M)$ such that
$$ \varphi (u,v,w)=g(u \times v,w)$$ for tangent vectors $u,v,w\in T(M)$.
For more information on $G_2$ manifolds see \cite{Joyce1}.

Next, we define an asymptotically cylindrical $G_2$ manifold with
one end.

Let $(M_0,\vp_0,g_0)$ be a $G_2$ manifold with calibration 3-form
$\vp_0$ and the Riemannian metric $g_0$. In this case Definitions
\ref{co2def1} and \ref{co2def22} take the following form, as in
\cite{Salu0}.
\begin{dfn} A $G_2$-manifold $(M_0,\vp_0,g_0)$ is called
{\it cylindrical} if $M_0=X\t\R$ and $(\vp_0,g_0)$ is compatible
with this product structure, that is,
\begin{equation*}
\vp_0=\Re\Om+\om\w\d t \qquad\text{and}\qquad g_0=g_X+\d t^2,
\end{equation*}
where $X$ is a compact, connected Calabi--Yau 3-fold with
K\"{a}hler form $\om$, Riemannian metric $g_X$ and holomorphic
(3,0)-form $\Om$. \label{co2def5}
\end{dfn}
\begin{dfn} A connected, complete $G_2$-manifold $(M,\vp,g)$ is
called {\it asymptotically cylindrical with decay rate} $\be$,
$\be<0$, if there exists a cylindrical $G_2$-manifold
$(M_0,\vp_0,g_0)$ with $M_0=X\t\R$ as above, a compact subset
$K\subset M$, a real number $R$, and a diffeomorphism
$\Psi:X\t(R,\iy)\ra M\sm K$ such that $\Psi^*(\vp)=\vp_0+\d\xi$
for some smooth 2-form $\xi$ on $X\t(R,\iy)$ with
$\bmd{\na^k\xi}=O(e^{\be t})$ on $X\t(R,\iy)$ for all $k\ge 0$,
where $\na_0$ is the Levi-Civita connection of the cylindrical
metric $g_0$. \label{co2def6}
\end{dfn}

Using similar relations between the geometric structures we can
define cylindrical and asymptotically cylindrical Calabi-Yau
manifolds. The details of these definitions and the applications
will appear in a forthcoming paper, \cite{Salu1}.

\subsection{The Laplacian operator on asymptotically cylindrical manifolds }
\label{co31}

In this section we prove Proposition \ref{co3prop37} and study the
properties of the cokernel of asymptotically cylindrical Laplacian
operator. This proposition plays an important role in Section \ref
{co3sec3} in the proof of Theorem \ref{co1thm}.

First we recall the basic facts about the properties of the
Laplacian operator and its (co)kernel on asymptotically
cylindrical manifolds with decay rate $\be<0$. More details can be
found in Lockhart and McOwen \cite{Lock}, \cite{LoMc}. Also one
can find similar results with a different analytical approach in
Melrose, \cite{Mel1}, \cite{Mel2}.

\begin{dfn} Let $(M,g)$ be an asymptotically cylindrical Riemannian manifold
asymptotic to $X\t\R$ with decay rate $\be$. Let
$\Delta_0:C^\iy(X\t\R)\ra C^\iy(X\t\R)$ be the {\it cylindrical}
Laplacian operator defined on the Riemannian cylinder which is the
asymptotic model for the ends of $M$ and invariant under
translations in $\R$. Let also $\Delta:C^\iy(M)\ra C^\iy(M)$ be
the Laplacian operator on $M$. Suppose that $\Psi$ is defined as
in Definition \ref{co2def22}. Then since $g$ is an asymptotically
cylindrical metric asymptotic to $g_0$, $\Delta$ is asymptotic to
$\Delta_0$ and $\Psi_*(\Delta)=\Delta_0+O(e^{\be t})$ as $t\ra\iy$
for $\be<0$ in the sense that their coefficients (in local
coordinates on $X$ and $\R$) are exponentially close. We call
$\Delta$ the {\it asymptotically cylindrical} Laplacian operator.

Choose a smooth function $\rho:M\ra\R$ such that $\Psi^*(\rho)
\equiv t$ on $X\t(R,\iy)$. This prescribes $\rho$ on $M\sm K$, so
we only have to extend $\rho$ over the compact set $K$. For $p\geq
1$, $k\geq 0$ and $\al\in\R$ we define the {\it weighted Sobolev
space} $L^p_{k,\al}(M)$ to be the set of functions $f$ on $M$ that
are locally integrable and $k$ times weakly differentiable and for
which the norm
\begin{equation}
\nm{f}_{L^p_{k,\al}}=\Bigl(\sum_{j=0}^{k}\int_M
e^{-\al\rho}\bmd{\na^jf}^p\d V\Bigr)^{1/p} \label{co3eq3}
\end{equation}
is finite. Then $L^p_{k,\al}(M)$ is a Banach space. Since $\rho$
is uniquely determined except on the compact set $K$, different
choices of $\rho$ give the same space $L^p_{k,\al}(M)$, with
equivalent norms. \label{co3def1}
\end{dfn}

Let $L^p_{k,\al}(M)$ be the weighted Sobolev spaces for some $\al$
such that $\be<\al<0$. Then one can show that $\Delta$ extends to
bounded linear operators
\begin{equation}
\Delta^p_{k,\al}:L^p_{k,\al}(M)\longra L^p_{k-2,\al}(M)
\label{co3eq0}
\end{equation}
for all $p>1$, $k\ge 2$ and $\al\in\R$.

In Definition \ref{co3def3} we define a set
$\D_{\Delta_0}\subset\R$, which determines the value of $\al$ so
that an asymptotically cylindrical operator is Fredholm.

\begin{dfn} Let $\Delta$ and $\Delta_0$
be the Laplacian operators on $M$ and $X\t\R$. Extend $\Delta_0$
to the complexifications $\Delta_0:C^\iy(X\t\R)\ot\C\ra
C^\iy(X\t\R)\ot\C$. Then $\D_{\Delta_0}$ is the set of $\al\in\R$
such that for some $\de\in\R$ there exists a nonzero section $s\in
C^\iy(X\t\R)\ot\C$ invariant under translations in $\R$ such that
$\Delta_0(e^{(\al+i\de)t}s)=0$. \label{co3def3}
\end{dfn}

Lockhart and McOwen prove \cite[Th.~1.1]{LoMc} that \eq{co3eq0} is
Fredholm if and only if $\al$ does not lie in a discrete set
$\D_{\Delta_0}\subset\R$.

\begin{thm} Let $(M,g)$ be an asymptotically cylindrical
Riemannian manifold asymptotic to $(X\t\R,g_0)$, and
$\Delta:C^\iy(M)\ra C^\iy(M)$ be the asymptotically cylindrical
Laplacian operator on $M$ between functions on $M$, asymptotic to
the cylindrical Laplacian operator $\Delta_0: C^\iy(X\t\R)\ra
C^\iy(X\t\R)$ on $X\t\R$. Let $\D_{\Delta_0}$ be defined as above.

Then $\D_{\Delta_0}$ is a discrete subset of $\R$, and for $p>1$,
$k\ge 2$ and $\al\in\R$, the extension
$\Delta^p_{k,\al}:L^p_{k,\al} (M)\ra L^p_{k-2,\al}(M)$ is Fredholm
if and only if $\al\notin\D_{\Delta_0}$. \label{co3thm1}
\end{thm}

Also by an elliptic regularity result \cite[Th.~3.7.2]{Lock} and
the weighted Sobolev Embedding Theorem \cite[Th.~3.10]{Lock} we
get
\begin{thm} For $\al\notin\D_{\Delta_0}$ the kernel $\Ker(\Delta^p_{k,\al})$
is independent of $p,k$, and is a finite-dimensional vector space
of smooth functions on $M$.
\label{co3prop1}
\end{thm}

\begin{thm} Let $\Delta^*$ denote the adjoint operator of $\Delta$. Then for all $\al\notin\D_{\Delta_0}$, $p,q>1$
with $\frac{1}{p}+\frac{1}{q}=1$ and $k,m\ge 2$ there is a natural
isomorphism
\begin{equation}
\Coker(\Delta^p_{k,\al})\cong\Ker\bigl((\Delta^*)^q_{m,-\al}\bigr)^*.
\label{co3eq6}
\end{equation}
\label{co3prop2}
\end{thm}
We now prove the following proposition. A related result with similar calculations can be found in Lockhart and McOwen, \cite[Th.~7.4]{LoMc}.
\begin{prop}  Let $M$ be an orientable, asymptotically cylindrical Riemannian manifold
with $l$ cylindrical ends $X_1\t(R,\iy ), \ldots, X_l\t(R,\iy )$.
Let $\Delta^p_{k,\al}$ be the asymptotically cylindrical Laplacian
operator on $M$ for some small $\al$, $\al<0$,
$[\al,0)\cap\D_{\Delta_0}=\emptyset$. Then
\begin{equation}
\dim(\Coker(\Delta^p_{k,\al}))=\dim(\Ker\bigl((\Delta^*)^q_{m,-\al}\bigr)^*)=l.
\end{equation}
\label{co3prop37}
\end{prop}
\begin{proof} When $\al\notin\D_{\Delta_0}$ we see from \eq{co3eq6} that

\begin{equation}
\begin{split}
\ind(\Delta^p_{k,\al})&=\dim\Ker(\Delta^p_{k,\al})
-\dim\Ker\bigl((\Delta^*)^q_{m,-\al}\bigr).\\
\ind(\Delta^q_{m,-\al})&=\dim\Ker(\Delta^q_{m,-\al})
-\dim\Ker\bigl((\Delta^*)^p_{k,\al}\bigr).
\end{split}
\label{co3eq7}
\end{equation}

Now since $\Delta$ is self adjoint, $\Delta=\Delta^*$, we have
\begin{equation}
\ind\bigl(\Delta^q_{m,-\al}\bigr)=
-\ind\bigl(\Delta^p_{k,\al}\bigr)=
\dim\Ker\bigl(\Delta^q_{m,-\al}\bigr)
-\dim\Ker\bigl(\Delta^p_{k,\al}\bigr) \label{co3eq78}
\end{equation}

Let $M$ be an asymptotically cylindrical Riemannian manifold with
one cylindrical end $X\t(R,\iy )$. Lockhart and McOwen show
\cite[Th.~6.2]{LoMc} that for $\al,\de\in\R\sm\D_{\Delta_0}$ with
$\al\le\de$ we have
\begin{equation}
\ind(\Delta^p_{k,\de})-\ind(\Delta^p_{k,\al})=
\sum_{\ep\in\D_{\Delta_0}:\al<\ep<\de}d(\ep), \label{co3eq8}
\end{equation}
where $d(\ep)\ge 1$ is the dimension of the vector space of
solutions $f\in C^\iy(X\t\R)\ot\C$ of a prescribed form as in Definition \ref{co3def3} with~$\Delta_0f=0$.

As $[\al,-\al]\cap\D_{\Delta_0}=\{0\}$, we can show that

\begin{equation}
\ind\bigl(\Delta^p_{k,-\al}\bigr)-
\ind\bigl(\Delta^p_{k,\al}\bigr) =\ts 2b^0(X). \label{co3eq141}
\end{equation}
By \eq{co3eq8}, $\ind\bigl(\Delta^q_{m,-\al}\bigr)-
\ind\bigl(\Delta^p_{k,\al}\bigr)$ is the dimension of the solution
space of $\Delta_0f=0$ on $X\t\R$ for $f$ polynomial in $t\in\R$.
By Hodge theory we deduce \eq{co3eq141}. By Theorems
\ref{co3prop1}, \ref{co3prop2} we can conclude that the dimensions
of kernel and cokernel are independent of $p,k$ and we can rewrite
(\ref{co3eq141}) as
\begin{equation}
\ind\bigl(\Delta^q_{m,-\al}\bigr)-
\ind\bigl(\Delta^p_{k,\al}\bigr) =\ts 2b^0(X). \label{co3eq15}
\end{equation}

Now let $M$ be an asymptotically cylindrical Riemannian manifold
with $l$ cylindrical ends $X_1\t(R,\iy ), \ldots, X_l\t(R,\iy )$.
In this case we get a contribution of 2 dimensions for each end
and the equation (\ref{co3eq15}) becomes
\begin{equation}
\ind\bigl(\Delta^q_{m,-\al}\bigr)-
\ind\bigl(\Delta^p_{k,\al}\bigr) =\ts 2l.
\label{co3eq151}
\end{equation}

Now since we have from equation (\ref{co3eq78}),
$\ind\bigl(\Delta^q_{m,-\al}\bigr)=
-\ind\bigl(\Delta^p_{k,\al}\bigr)$ we get
\begin{equation}
\ind\bigl(\Delta^q_{m,-\al}\bigr)
=\dim\Ker\bigl(\Delta^q_{m,-\al}\bigr)-\dim\Ker\bigl(\Delta^q_{m,\al}\bigr)=l.
\label{co3eq152}
\end{equation}

But by maximum principle the harmonic functions which decay with
rate $\al$, $\al<0$ cannot be nonzero functions and one can easily
see that $\Ker(\Delta^q_{m,\al})=0$. This with equation
(\ref{co3eq152}) yields the result.
\end{proof}

\section{Proof of Theorem \ref{co1thm}}
\label{co3sec3}
%\subsection{Construction of the constant harmonic 1-form $\ga$}
In order to prove Theorem \ref{co1thm} we need to construct a
closed and co-closed one form $\ga=\d f$ for some harmonic function
$f$ on $M$ and show that $\nabla\ga=0$.

\subsection{Construction of $f$ on $M$.}
We start with constructing a function $f$ which is harmonic on $M$
with $l$ cylindrical ends $X_1\t(R,\iy ), \ldots, X_l\t(R,\iy )$.
Let $f_0$ be a smooth function on $M$ of the form $f_0=C_it_i+D_i$
on $X_i\t(R,\iy )$ for each $i=1,\ldots,l$. Here $t_1,\ldots, t_l$
are coordinates on cylindrical ends and $C_i,D_i\in\R$. The
function $f_0$ smoothly interpolates between the cylindrical ends.
The harmonic function $f$ on $M$ is of the form $f=f_0+f'$ for
some function $f'$. Using elliptic estimates, \cite[Th.~3.7.2]{Lock} and embeddings of Sobolev spaces into H\"{o}lder spaces \cite[Sec.2.7.]{Aubin} one can take $f'\in L^p_{k,\al}(M)$ to be smooth which decays like $O(e^{\al t})$. Therefore $f=C_it_i+D_i+O(e^{\al t_i})$, $\al<0$,
on each end $X_i\t\R$.

In order to show the existence of the function $f$ such that
$f=f_0+f'$ we need to find a function $f'$ that satisfies
\begin{equation}
\Delta f'=-\Delta f_0.
\end{equation}
This is equivalent to showing that $\Delta f_0\bot
\Coker(\Delta^p_{k,\al})$ and in Proposition \ref{co3prop41} we
prove this fact. Let $f_0=f_0^{(C_i,D_i)}$ be defined as above.

\begin{prop} Let $M$ be an asymptotically cylindrical Riemannian manifold
with $l$ cylindrical ends $X_1\t(R,\iy ), \ldots, X_l\t(R,\iy )$.
Also for $h\in (\Coker(\Delta^p_{k,\al}))^*$ let $\Phi$ be a map
defined as
\begin{equation}
\begin{split}
\Phi:(C_i,D_i)&\longrightarrow \Coker(\Delta^p_{k,\al})\\
\Phi(C_i,D_i):h&\longmapsto\langle\Delta
f_0^{(C_i,D_i)},h\rangle_{L^2}
\end{split}
\end{equation}
\noindent where $\langle\cdot,\cdot\rangle_{L^2}$ is the
$L^2$-inner product defined on the space of functions on $M$.

Then for any element $(C_i,D_i)\in\Ker\Phi$ the map
$\Ker\Phi\longrightarrow(\Coker(\Delta^p_{k,\al}))^*$ given by
\begin{equation}
(C_i,D_i)\longrightarrow \langle f^{(C_i,D_i)}| \Delta
f^{(C_i,D_i)}=0, f^{(C_i,D_i)}=C_i t_i+D_i+O(e^{\al t_i}) {\rm
\hspace{.1in} on \hspace{.1in}} X_i\t(R,\iy)\rangle
\end{equation}
is an isomorphism. \label{co3prop41}
\end{prop}
\begin{proof}
First note that $\Phi$ is linear in $(C_i,D_i)$ as we defined
$f_0^{(C_i,D_i)}$ to be linear in $(C_i,D_i)$. Since there are
2$l$ parameters in $(C_i,D_i)$ and also
$\dim\Coker(\Delta^p_{k,\al})=l$ from Proposition \ref{co3prop37},
$\Phi$ is a linear map between $\R^{2l}\rightarrow \R^l$ which
guarantees that $\dim\Ker\Phi\geq l$. But for each
$(C_i,D_i)\in\Ker\Phi$, if $\Delta f_0^{(C_i,D_i)}\bot
\Coker(\Delta^p_{k,\al})$ then $(f')^{(C_i,D_i)}\in
L^p_{k,\al}(M)$ exists and this implies that there exists an
$f^{(C_i,D_i)}$ such that $\Delta f^{(C_i,D_i)}=0$ and
$f^{(C_i,D_i)}=C_i t_i+D_i+O(e^{\al t_i})$ on i'th end.

Therefore we have an injective linear map
\begin{equation*}
\Ker\Phi\longrightarrow\langle f^{(C_i,D_i)}| \Delta
f^{(C_i,D_i)}=0, f^{(C_i,D_i)}=C_i t_i+D_i+O(e^{\al
t_i})\rangle\subset(\Coker(\Delta^p_{k,\al}))^*
\end{equation*}

\noindent and hence we get an injective linear map
\begin{equation*}
\begin{split}
\Ker\Phi &\longrightarrow (\Coker(\Delta^p_{k,\al}))^*\\
(C_i,D_i) &\longmapsto
f^{(C_i,D_i)}=f_0^{(C_i,D_i)}+(f')^{(C_i,D_i)}.
\end{split}
\end{equation*}

This implies that $\dim\Ker\Phi\leq \dim (\Coker(\Delta^p_{k,\al}))^*
=l$ and hence $\dim\Ker\Phi=l$. Therefore $\Ker\Phi\rightarrow
(\Coker(\Delta^p_{k,\al}))^*$ is an isomorphism as an injective linear
map between spaces of dimension $l$ and this completes the proof.
\end{proof}

\subsection{Existence of constant harmonic 1-form $\ga=df$ on $M$.}
\label{co3sec2}
 After constructing the harmonic function $f$, one
can easily obtain a closed and coclosed 1-form $\ga=df$ on $M$. We
now show that $\ga=df$ is constant for Ricci-flat manifolds. Let
$g$ be the Riemannian metric on $M$. In the index notation the
{\em Weitzenb\"{o}ck formula} for a 1-form $\xi$ is
\begin{equation}
(\d\d^*+\d^*\d)\xi_a=\nabla^*\nabla\xi_a+R_{ab}g^{bc}\xi_c.
\label{co3eq222}
\end{equation}
where $R_{ab}$ is the Ricci curvature of $g$. Suppose that
$R_{ab}=0$ then we have
\begin{equation}
(\d\d^*+\d^*\d)\xi=\nabla^*\nabla\xi \label{co3eq2}
\end{equation}
on 1-forms on a Ricci-flat manifold $M$.
\begin{prop} Let $M$ be an asymptotically cylindrical Riemannian manifold with Ricci-flat
metric $g$. Let $\ga=df$ be a harmonic 1-form defined as above.
Then $\ga=df$ is constant, that is $\nabla\ga=0$.
\label{co3prop88}
\end{prop}
\begin{proof} Let $K_R= M\backslash((X_1\t(R,\iy )\cup\ldots\cup (X_l\t(R,\iy ))$,
$R\gg0$. Then $\partial K_R=(X_1\t \{R\})\cup\ldots\cup (X_l\t
\{R\})$. Since $\ga$ is harmonic we have $(\d\d^*+\d^*\d)\ga=0$
and by equation (\ref{co3eq2}) we get
\begin{equation}
\int_{K_R}\ga((\d\d^*+\d^*\d)\ga)=\int_{K_R}\ga(\nabla^*\nabla\ga)=0
\end{equation}
\noindent and by Stoke's Theorem,
\begin{equation}
\int_{K_R}|\nabla\ga|^2=\int_{X_1\t \{R\}}\ga\cdot(\nabla\ga\cdot
\eta_1)+\ldots+\int_{X_l\t \{R\}}\ga\cdot(\nabla\ga\cdot \eta_l).
\label{co3eq455}
\end{equation}
for normal vectors $\eta_1,\ldots,\eta_l$. The asymptotic
behaviour of $\ga=df$ guarantees that the right hand side of
equation (\ref{co3eq455}) goes to zero as $R\rightarrow \infty$
and yields the result.
\end{proof}

We now complete the proof of Theorem \ref{co1thm}. Let $(M,g)$ be
an asymptotically cylindrical Ricci-flat manifold with $l$ ends.
Then $(M,g)$ is asymptotic to the cylinders $X_1\t(R,\iy )$,...,
$X_l\t(R,\iy )$, for $R\gg0$. Proposition \ref{co3prop41}
guarantees the existence of an harmonic function $f$ on $M$.

Note that if $l=1$, by Proposition \ref{co3prop37},
dim($\Coker(\Delta^p_{k,\al}))=1$ and $f$ should be constant and
so $df=0$. If $l\ge 2$ then one can choose $f$ not constant and
construct a nonzero, closed and coclosed 1-form $\ga=df$. Since
$(M,g)$ is Ricci-flat, by Proposition \ref{co3prop88}, $\ga$ is
constant. Thus the results \cite[Sec.3.2]{Joyce1} show that the
metric $g$ on $M$ is locally reducible. Rescale $f$ if necessary
so that $\vert\gamma\vert\equiv 1$.

Also note that by Proposition \ref{co3prop88}, $df$ is covariant
constant and so $|df|$ is constant. Therefore $\ga=df$ never
vanishes on $M$, and $f$ can be viewed as a Morse function with no
critical points. Define $X$ to be the level set $f^{-1}(0)$. This
is a nonsingular submanifold of M as $df\ne 0$ everywhere, and is
compact as $f\rightarrow\pm\infty$ at the ends of $M$. This
implies that $M$ is at least a topological product $M=X\t\R$.

We now construct a diffeomorphism $X\times\R\rightarrow M$. For
$t\in\R$ define $\upsilon_t:M\rightarrow M$ to be the gradient
flow of $f$ for time $t$. As $M$ is complete and  $\vert
df\vert\equiv 1$ this is a well-defined diffeomorphism, which
induces a diffeomorphism between $X$ and the level set
$f^{-1}(t)$. Define $\Upsilon:X\times\R\rightarrow M$ by
$\Upsilon(x,t)=\upsilon_t(x)$. This is a bijection as for each
$t$, $x\mapsto\Upsilon(x,t)=\upsilon_t(x)$ is a bijection between
$X$ and $f^{-1}(t)$, and it is easy to see $\Upsilon$ is a
diffeomorphism.

As $\vert df\vert\equiv 1$, the metric $\Upsilon^*(g)$ is of the
form $g(x,t)+dt^2$, where $g(x,t)$ is a 1-parameter family of
metrics on $X$ depending on $t\in\R$. Furthermore, since
$\Upsilon^*(g)$ is locally reducible and $dt$ is constant, we see
that $g(x,t)$ is locally and hence globally independent of $t$.
Thus $M$ is a Riemannian product $X\times\R$. Since $M$ is
connected $X$ is connected, so $M$ has exactly two ends, so $l=2$
and it is a cylinder.

\begin{rem} If $(M,g)$ is an asymptotically cylindrical $G_2$
manifold with two cylindrical ends, then Hol$(g)$ preserves $\ga$
and there is a reduction in the holonomy group and hence Hol$(g)\subset SU(3)$. This implies that
$M=X\times\R$ where $X$ is a Calabi-Yau 3-fold and hence $M$ is a
cylinder.

\end{rem}

{\small{\it Acknowledgements.} The author is grateful to Dominic
Joyce for his help and encouragement during this project and to
the AWM for their grant support. Special thanks to Jeff Viaclovsky
and an anonymous referee for many valuable comments.
 }

\end{document}